\documentclass[12pt,a4paper]{article}

\usepackage{latexsym}
\usepackage{amsmath}
\usepackage{amssymb}
\usepackage{arydshln}

\usepackage{cite}
\usepackage{stmaryrd}
\usepackage{enumerate}

\usepackage{hyperref}

\usepackage{color}
\usepackage{lineno}
\usepackage{graphicx}
\usepackage{ae}
\usepackage{amsmath}
\usepackage{amssymb}
\usepackage{latexsym}
\usepackage{url}
\usepackage{epsfig}
\usepackage{mathrsfs}
\usepackage{amsfonts}
\usepackage{amsthm}

\usepackage{float}
\usepackage{subfig}
\usepackage{tikz}
\usetikzlibrary{positioning}

\newtheorem{theorem}{Theorem}[section]

\newtheorem{lemma}[theorem]{Lemma}

\theoremstyle{definition}

\numberwithin{equation}{section} 
\allowdisplaybreaks    

\long\def\delete#1{}

\usetikzlibrary{decorations.markings}

\tikzstyle{vertex}=[circle, draw, inner sep=0pt, minimum size=6pt]
\tikzstyle{directed}=[postaction={decorate,
	decoration={markings,mark=at position 0.3 with {\arrow{stealth}}}
}]

\usepackage{xcolor}
\usepackage[normalem]{ulem}

\begin{document}
\title {Index perturbation of signed graphs}

\author{Linfeng Xie$^{a,b}$,~Xiaogang Liu$^{a,b,}$\thanks{Supported by the National Natural Science Foundation of China (No. 12371358).}~$^,$\thanks{ Corresponding author. Email addresses: xielinfeng@mail.nwpu.edu.cn, xiaogliu@nwpu.edu.cn}~
	\\[2mm]
	{\small $^a$School of Mathematics and Statistics,}\\[-0.8ex]
	{\small Northwestern Polytechnical University, Xi'an, Shaanxi 710072, P.R.~China}\\
	{\small $^b$Xi'an-Budapest Joint Research Center for Combinatorics,}\\[-0.8ex]
	{\small Northwestern Polytechnical University, Xi'an, Shaanxi 710129, P.R. China}\\
}
\date{}

\openup 0.5\jot
\maketitle

\begin{abstract}
	Let $\Gamma = (G, \sigma)$ be a signed graph and $v$ a non-isolated vertex of $\Gamma$. Let $\Gamma-v$ denote the graph obtained by deleting the vertex $v$ together with all signed edges incident to it from $\Gamma$, and $d_{\Gamma}(v)$ the degree of $v$ in $\Gamma$. In this paper, we prove that the largest eigenvalue $\lambda_1(\Gamma)$ of $\Gamma$ satisfies
	\[
	\lambda_1(\Gamma) \le \sqrt{\lambda_1^2(\Gamma - v) + 2d_\Gamma(v) - 1},
	\]
	and we also present a refined version of this bound. Moreover, we characterize the extremal signed graphs achieving equality when $\Gamma$ is connected and $d_\Gamma(v)\ge 2$, which are switching equivalent to the balanced complete signed graph. 

	\smallskip
	
	
\end{abstract}

\section{Introduction}
A \emph{signed graph} $\Gamma=(G,\sigma)$ is a pair $(G,\sigma)$, where $G=(V(G),E(G))$ is a simple graph and $\sigma : E  \to \left \{+1, -1 \right \}$ is the sign function defined on the edge set $E = E(G)$. The unsigned graph $G$ is called the \emph{underlying graph} of $\Gamma$. In 1946, Heider first introduced the concept of signed graphs in his study of balance theory in social psychology \cite{Hei}. Heider sought to explain how friendly/hostile relationships in social networks influence group stability. Later, in 1953, Harary formalized this theory and established the mathematical foundation of signed graphs \cite{Har}. For more information of signed graphs, please refer to \cite{Cam, Zaslavsky1982, Zaslavsky2008,Zaslavsky2018}.

Let $\Gamma=(G,\sigma)$ be a signed graph with the vertex set $V(\Gamma)=\left \{v_{1},v_{2},\dots ,v_{n}\right \}$ and the edge set $E(\Gamma)=\left \{e_{1},e_{2},\dots ,e_{m}\right \}$. Denote by $|V(\Gamma)|$ and $|E(\Gamma)|$ the order and the size of $\Gamma$, respectively. Let $N_{\Gamma}(v_i)$ denote the set of neighbors of a vertex $v_i$ in $\Gamma$ and $d_{\Gamma}(v_i)=|N_{\Gamma}(v_i)|$ the degree of $v_i$ in $\Gamma$. An edge $v_i v_j$ is called a \emph{positive edge} (respectively, \emph{negative edge}) if $\sigma(v_i v_j)=+1$ (respectively, $\sigma(v_i v_j)=-1$). The set of all positive edges (respectively, negative edges) in $\Gamma$ is denoted by $E^{+}(\Gamma)$ (respectively, $E^{-}(\Gamma)$). Denote by $(G, +)$ (respectively, $(G, -)$) the signed graph whose edges are all positive (respectively, negative). The \emph{adjacency matrix} of $\Gamma$ is defined as $A(\Gamma)=(a_{ij}^{\sigma})$, where $a_{ij}^{\sigma}=\sigma(v_{i}v_{j})a_{ij}$ and $a_{ij}=1$ if $ v_{i} $ and $v_{j}$ are adjacent, and $a_{ij}=0$ otherwise. The eigenvalues of $A(\Gamma)$ are denoted by
$$
\lambda_{1}(\Gamma)\ge \lambda_{2}(\Gamma)\ge \dots \ge \lambda_{n}(\Gamma),
$$
which are called the \emph{adjacency spectrum} of $\Gamma$. The \emph{index} of $\Gamma$ is $\lambda_{1}(\Gamma)$.

Let $\emptyset \ne U\subseteq V(G)$. The operation of changing the signs of all edges between $U$ and $V(G) \setminus U$ is called a \emph{switching}, and is also referred to as \emph{switching $\Gamma$ at a vertex subset $U$}. A signed graph $\Gamma^{\prime}$ is said to be \emph{switching equivalent} to $\Gamma$ if $ \Gamma^{\prime}$ is obtained from $\Gamma$ by a finite sequence of switchings, denoted as $\Gamma\sim \Gamma^{\prime}$. A cycle is called \emph{positive} if the number of its negative edges is even; otherwise, \emph{negative}. A signed graph is called \emph{balanced} if each of its cycles is positive; otherwise, \emph{unbalanced}.

For a graph $G$, let $\lambda_{1}(G)$ be the largest eigenvalue of $G$. For any $v\in V(G)$, let $G-v$ denote the graph obtained by deleting the vertex $v$ together with all edges incident to it from $G$. By \cite[Proposition~1.3.9]{Cvt},
$$
\lambda_{1}(G-v)<\lambda_{1}(G)
$$
holds for all connected graphs $G$ and any vertex $v\in V(G)$. Hence, investigating the difference $\lambda_{1}(G)-\lambda_{1}(G-v)$ is an interesting problem in spectral graph theory. In 2010, Nikiforov \cite{Nikiforov-LAA-2010} proved that
\[
\lambda_1(G) \leq \lambda_1(G - v_k) \frac{1 - x_k^2}{1 - 2x_k^2},
\]
where $x_k$ is the smallest entry of the Perron unit eigenvector of $G$ and $v_k$ is the vertex with respect to $x_k$. Later, in 2019, Guo, Wang and Li \cite{Guo-Wang-Li-DM-2019} proposed a conjecture that for any graph $G$ and a non-isolated vertex $v \in V(G)$,
$$
\lambda_1(G) \le \sqrt{\lambda^2_1(G-v) + 2d_{G}(v) - 1},
$$
and verified the inequality for the case $d_{G}(v)=1$, where $d_{G}(v)$ denotes the degree of $v$ in $G$. In the following year, Sun and Das \cite{Sun-Das-LAA-2020} confirmed the aforementioned conjecture, and they characterized all connected graphs for which this bound is attained (See Theorem \ref{TDeleteV-11}). Recently, Liu and Ning \cite{Liu-Ning-AR-2026.3.31} gave an alternative short proof of the result.
\begin{theorem}\emph{(See \cite{Sun-Das-LAA-2020,Liu-Ning-AR-2026.3.31})}\label{TDeleteV-11}
	Let $G$ be a simple graph and $v \in V(G)$ a non-isolated vertex. Then
	\[
	\lambda_1(G) \le \sqrt{\lambda^2_1(G-v) + 2d_{G}(v) - 1}.
	\]
	If $G$ is connected, then equality holds if and only if $G \cong K_{1,n-1}$ with $d_{G}(v)=1$, or $G \cong K_n$.
	
\end{theorem}

For a signed graph $\Gamma$ and any $v\in V(\Gamma)$, let $\Gamma-v$ denote the graph obtained by deleting the vertex $v$ together with all signed edges incident to it from $\Gamma$. Let $\epsilon(\Gamma)$ denote the \emph{frustration index} of $\Gamma$, which is the minimum number of edges to remove for balance.

In this paper, we extend Theorem \ref{TDeleteV-11} to signed graphs, which is stated as follows. The method used to prove our main result is inspired by the work in \cite{Liu-Ning-AR-2026.3.31}.

\begin{theorem}\label{largest eigenvalue perturbation}
	Let $\Gamma=(G,\sigma)$ be a signed graph and $v\in V(\Gamma)$ a non-isolated vertex. Then
	\begin{equation}\label{largest eigenvalue perturbation-theorem-1}
		\lambda_{1}(\Gamma)\le \sqrt{\lambda^2_{1}(\Gamma-v)+2d_{\Gamma}(v)-1}.
	\end{equation}
	If $\Gamma$ is connected and $d_{\Gamma}(v)\ge 2$, then equality holds in \eqref{largest eigenvalue perturbation-theorem-1} if and only if $\Gamma\sim (K_{n},+)$.
	
	Furthermore, if $d_{\Gamma}(v)\ge 2$ and $\lambda_{1}(\Gamma-v)<\lambda_{1}(\Gamma)$, then
	\begin{equation}\label{largest eigenvalue perturbation-theorem-2}
		\lambda_{1}(\Gamma)\le \sqrt{\lambda^2_{1}(\Gamma-v)+2d_{\Gamma}(v)-1-\frac{8\epsilon(\Gamma)}{\lambda_{1}(\Gamma)(d_{\Gamma}(v)+1)}}.
	\end{equation}
	
\end{theorem}

\section{Proof of Theorem \ref{largest eigenvalue perturbation}}\label{largest eigenvalue perturbation-1}
Let $\Gamma=(G,\sigma)$ be a signed graph. For $U\subseteq V(\Gamma)$, let $\Gamma[U]$ denote the signed subgraph of $\Gamma$ induced by $U$, with edge signs inherited from $\Gamma$. Sometimes, we say that $U$ induces $\Gamma[U]$. Let $A$ and $B$ be two symmetric matrices. We write $A \succeq B$ if $A-B$ is positive semi-definite.

\begin{lemma}\emph{(See \cite{Cvet})}\label{Cvet}
	Let $\Gamma$ be a signed graph and $U$ a subset of $V(\Gamma)$ with $|U|=k$. Then
	$$
	\lambda_{i}(\Gamma) \ge \lambda_{i}(\Gamma[V(\Gamma)\setminus U]) \ge \lambda_{i+k}(\Gamma),~~\text{for $1 \le i \le n-k$}.
	$$
\end{lemma}

\begin{Tproof}\textbf{~of~Theorem~\ref{largest eigenvalue perturbation}.}
	For brevity, set $\lambda_{1}:= \lambda_1(\Gamma)$ and $\mu_1 := \lambda_1(\Gamma-v)$. By Lemma \ref{Cvet}, $\lambda_1 \ge \mu_1$. If $\lambda_1 = \mu_1$, then \eqref{largest eigenvalue perturbation-theorem-1} follows immediately. In the following, we assume that $\lambda_1 >\mu_1$.
	
	By switching equivalence, there exists a signed graph $\Gamma'$ such hat $\Gamma'\sim \Gamma$ and $A(\Gamma')$ has the block form
	\[
	A(\Gamma') = \begin{bmatrix}
		0 & \mathbf{a} ^T \\
		\mathbf{a} & B
	\end{bmatrix},
	\]
	where $B = A(\Gamma'-v)$ and $\mathbf{a}$ is the $0$-$1$ column vector indicating the neighbors of $v$ in $\Gamma'-v$. Clearly, $\mathbf{a}^T\mathbf{a} = d_{\Gamma}(v)$.
	
	Recall that $\mu_1<\lambda_1$. We obtain that $\lambda_1 I - B$ is positive definite and invertible. Hence, $\det(\lambda_1 I - B) \ne 0$. By a simple calculation, we have
	\[
	\left(\lambda_1 - \mathbf{a}^T(\lambda_1 I - B)^{-1}\mathbf{a}\right)\det(\lambda_1 I - B) =0.
	\]
	So, $\lambda_1 = \mathbf{a}^T(\lambda_1 I - B)^{-1}\mathbf{a}$. Let $t$ be an arbitrary eigenvalue of $B$. Since $|t| \le \mu_1 < \lambda_1$, we have
	\[
	\frac{\lambda_1 + t}{\lambda^2_1 - \mu^2_1} - \frac{1}{\lambda_1 - t}
	= \frac{\mu^2_1 - t^2}{(\lambda^2_1 -\mu^2_1)(\lambda_1 - t)} \ge 0.
	\]
	Then
	$$
	\frac{\lambda_1 I + B}{\lambda^2_1 - \mu^2_1} \succeq (\lambda_1 I - B)^{-1}.
	$$
	Recall that  $\lambda_1 I - B$ is positive definite. Then there exists an orthogonal matrix $U$ such that $B=U\operatorname{diag}\left \{t_1,t_2,\dots,t_{n-1}\right \}U^T$. Hence,
	\begin{align*}
		M :=& \frac{\lambda_1 I + B}{\lambda^2_1 - \mu^2_1} - (\lambda_1 I - B)^{-1}\\
		=& U\operatorname{diag}\left \{\frac{\lambda_1 + t_1}{\lambda^2_1 - \mu^2_1} - \frac{1}{\lambda_1 - t_1},
		\dots,
		\frac{\lambda_1 + t_{n-1}}{\lambda^2_1 - \mu^2_1} - \frac{1}{\lambda_1 - t_{n-1}}  \right \}U^T \\
		\succeq &\mathbf{0}.
	\end{align*}
	Recall that $\lambda_1 = \mathbf{a}^T(\lambda_1 I - B)^{-1}\mathbf{a}$. Then
	\begin{equation}\label{perturbation-equation}
		\frac{\lambda_1 d_{\Gamma}(v) + \mathbf{a}^T B\mathbf{a}}{\lambda^2_1 - \mu^2_1} - \lambda_1
		= \frac{\lambda_1 \mathbf{a}^T\mathbf{a} + \mathbf{a}^T B\mathbf{a}}{\lambda^2_1 - \mu^2_1}
		- \mathbf{a}^T (\lambda_1 I - B)^{-1}\mathbf{a}
		= \mathbf{a}^T M\mathbf{a} \ge 0.
	\end{equation}

	Set $m^+:=|E^+(\Gamma'[N_{\Gamma'}(v)])|$ and $m^-:=|E^-(\Gamma'[N_{\Gamma'}(v)])|$. Clearly, $\mathbf{a}^T B\mathbf{a}=2(m^+-m^-)$. If $d_{\Gamma}(v)=1$, then $\mathbf{a}^T B\mathbf{a}=0$. By \eqref{perturbation-equation}, we obtain that
	$$
	\lambda^2_1 - \mu^2_1\le d_{\Gamma}(v)=1=2d_{\Gamma}(v)-1,
	$$
	which leads to \eqref{largest eigenvalue perturbation-theorem-1}.
	
	In the following, we consider the case $d_{\Gamma}(v)\ge 2$.
	
	Let $\Sigma=\Gamma[N_{\Gamma'}(v)\cup \left\{v\right\}]$. Then, $|V(\Sigma)|=d_{\Gamma}(v)+1$ and $2(|E^+(\Sigma)|-|E^-(\Sigma)|)=2(d_{\Gamma}(v)+m^+-m^-)$. Hence, by Rayleigh principle,
	$$
	\lambda_{1}(\Sigma)\ge \frac{2(d_{\Gamma}(v)+m^+-m^-)}{d_{\Gamma}(v)+1}.
	$$
	By a simple calculation, we have
	\begin{align*}
		\frac{2(d_{\Gamma}(v)+m^+-m^-)}{d_{\Gamma}(v)+1}-\frac{2(m^+-m^-)}{d_{\Gamma}(v)-1} &= \frac{2\left(d_{\Gamma}(v)(d_{\Gamma}(v)-1)-2(m^+-m^-)\right)}{d_{\Gamma}^2(v)-1}\\
		&= \frac{2\left(d_{\Gamma}(v)(d_{\Gamma}(v)-1)-2(m^++m^-)+4m^-\right)}{d_{\Gamma}^2(v)-1}.
	\end{align*}
	Note that $2(m^++m^-)\le d_{\Gamma}(v)(d_{\Gamma}(v)-1)$ and $m^-\ge \epsilon(\Gamma)$. Then
	$$
	\frac{2(d_{\Gamma}(v)+m^+-m^-)}{d_{\Gamma}(v)+1}-\frac{2(m^+-m^-)}{d_{\Gamma}(v)-1}\ge \frac{8\epsilon(\Gamma)}{d_{\Gamma}^2(v)-1}.
	$$
	By Lemma \ref{Cvet},
	$$
	\lambda_{1}\ge \lambda_{1}(\Sigma)\ge \frac{2(m^+-m^-)}{d_{\Gamma}(v)-1}+\frac{8\epsilon(\Gamma)}{d_{\Gamma}^2(v)-1}.
	$$
	Hence,
	$$
	\mathbf{a}^T B\mathbf{a}=2(m^+-m^-)\le \lambda_{1}\left(d_{\Gamma}(v)-1\right)-\frac{8\epsilon(\Gamma)}{d_{\Gamma}(v)+1}.
	$$
	By \eqref{perturbation-equation}, we obtain that
	$$
	\lambda^2_1-\mu^2_1\le 2d_{\Gamma}(v)-1-\frac{8\epsilon(\Gamma)}{\lambda_{1}(\Gamma)\left(d_{\Gamma}(v)+1\right)}\le 2d_{\Gamma}(v)-1.
	$$
	This yields \eqref{largest eigenvalue perturbation-theorem-2}.
	
	Finally, we characterize the graphs attaining the upper bounds in \eqref{largest eigenvalue perturbation-theorem-1} when $\Gamma'$ is connected and $d_{\Gamma}(v)\ge 2$. From the above proof, we have (i) $\lambda_{1}=\lambda_{1}(\Sigma)$, (ii) $2(m^++m^-)= d_{\Gamma}(v)(d_{\Gamma}(v)-1)$ and (iii) $2d_{\Gamma}(v)-1-\frac{8\epsilon(\Gamma)}{\lambda_{1}(\Gamma)\left(d_{\Gamma}(v)+1\right)}= 2d_{\Gamma}(v)-1$. By (iii), we have $\epsilon(\Gamma)=0$ and hence $\Gamma'$ is balanced. Consequently, by (i) and connectivity, $\Gamma'\cong \Sigma$. By (ii), we obtain that $\Sigma$ is a signed graph whose underlying graph is complete graph. Thus, $\Gamma$ is switching equivalent to a balanced complete signed graph.
\end{Tproof}


\begin{thebibliography}{99}
	\small{
		\bibitem{Cam}
		P. J. Cameron, J. J. Seidel, S. V. Tsaranov, Signed graphs, root lattices, and Coxeter groups, J. Algebra 164 (1) (1994) 173--209.	
		
		\bibitem{Cvet}
		D. M. Cvetkovi{\'c}, M. Doob, H. Sachs, Spectra of Graphs, Academic Press, New York, 1980.	
		
		\bibitem{Cvt}
		D. M. Cvetkovi{\'c}, P. Rowlinson, S. Simi{\'c}, An Introduction to the Theory of Graph Spectra, Cambridge University Press, Cambridge, 2010.
		
		\bibitem{Har}
		F. Harary, On the notion of balance of a signed graph, Michigan Math. J. 2 (2) (1953) 143--146.
		
		\bibitem{Hei}
		F. Heider, Attitudes and cognitive organization, J. Psychol. 21 (1) (1946) 107--112.
		
		\bibitem{Guo-Wang-Li-DM-2019}
		J.-M. Guo, Z.-W. Wang, X. Li, Sharp upper bounds of the spectral radius of a graph, Discrete Math. 342 (2019) 2559--2563.
		
		\bibitem{Liu-Ning-AR-2026.3.31}
		L. Liu, B. Ning, A short proof of an inequality on spectral radius perturbation, arXiv:2603.29335 (2026).
		
		\bibitem{Sun-Das-LAA-2020}
		S. Sun, K. C. Das, A conjecture on the spectral radius of graphs, Linear Algebra Appl. 588 (2020) 74--80.
		
		\bibitem{Nikiforov-LAA-2010}
		V. Nikiforov, The spectral radius of graphs without paths and cycles of specified length,  Linear Algebra Appl. 432 (9) (2010) 2243--2256.
		
		\bibitem{Zaslavsky1982}
		T. Zaslavsky, Signed graphs, Discrete Appl. Math. 4 (1) (1982) 47--74.
		
		\bibitem{Zaslavsky2008}
		T. Zaslavsky, Matrices in the theory of signed simple graphs, in Advances in Discrete Mathematics and Applications: Mysore, 2008, 207--229, Ramanujan Math. Soc. Lect. Notes Ser., 13, Ramanujan Math. Soc., Mysore, 2010.
		
		\bibitem{Zaslavsky2018}
		T. Zaslavsky, A mathematical bibliography of signed and gain graphs and allied areas, Electron. J. Combin. 5 (2018), Dynamic Surveys 8, 524 pages.
		
	}
\end{thebibliography}
\end{document}